\documentclass[a4paper]{article}

\usepackage{geometry}
\usepackage[english]{babel}
\usepackage{epsfig}
\usepackage{amsmath}
\usepackage{amssymb}
\usepackage{amsthm}

\usepackage{amsthm}
\newtheorem{theorem}{Theorem}[section]
\newtheorem{definition}[theorem]{Definition}
\newtheorem{lemma}[theorem]{Lemma}
\newtheorem{proposition}[theorem]{Proposition}
\newtheorem{corollary}[theorem]{Corollary}
\newtheorem{remark}[theorem]{Remark}
\newtheorem{example}[theorem]{Example}

\usepackage{amsfonts}

\newcommand{\hh}{{\mathbb{H}}}

\newcommand{\cc}{{\mathbb{C}}}
\newcommand{\rr}{{\mathbb{R}}}
\newcommand{\zz}{{\mathbb{Z}}}
\newcommand{\nn}{{\mathbb{N}}}

\newcommand{\s}{{\mathbb{S}}}
\newcommand{\ext}{\textnormal{ext}}

\title{\bf Singularities of slice regular functions}
\author{Caterina Stoppato\footnote{Partially supported by GNSAGA of the INdAM, by PRIN ``Geometria Differenziale e Analisi Globale'' and by PRIN ``Propriet\`a geometriche delle variet\`a reali e complesse'' of the MIUR.} \\ 
\small Dipartimento di Matematica ``U. Dini'', Universit\`a di Firenze \\ 
\small Viale Morgagni 67/A, 50134 Firenze, Italy\\
\small stoppato@math.unifi.it\\}

\date{  }

\begin{document}
\maketitle


\begin{abstract}
Beginning in 2006, G. Gentili and D.C. Struppa developed a theory of regular quaternionic functions with properties that  recall classical results in complex analysis. For instance, in each Euclidean ball $B(0,R)$ centered at $0$ the set of regular functions coincides with that of quaternionic power series $\sum_{n \in \nn} q^n a_n$ converging $B(0,R)$. 
In 2009 the author proposed a classification of singularities of regular functions as removable, essential or as poles and studied poles by constructing the ring of quotients. In that article, not only the statements, but also the proving techniques were confined to the special case of balls $B(0,R)$. 
In a subsequent paper, F. Colombo, G. Gentili, I. Sabadini and D.C. Struppa (2009) identified a larger class of domains, on which the theory of regular functions is natural and not limited to quaternionic power series. 
The present article studies singularities in this new context, beginning with the construction of the ring of quotients and of Laurent-type expansions at points $p$ other than the origin. These expansions, which differ significantly from their complex analogs,  allow a classification of singularities that is consistent with the one given in 2009. Poles are studied, as well as essential singularities, for which a version of the Casorati-Weierstrass Theorem is proven.
\end{abstract}



\section{Introduction}

Let $\hh$ denote the real algebra of quaternions, and recall that it consists of the vector space $\rr^4$ on which multiplication is constructed as follows: for the standard basis of $\rr^4$, denoted $1,i,j,k$ one defines
$$i^2=j^2 = k^2 = -1,$$
$$ij=-ji =k, jk=-kj=i, ki=-ik=j,$$
and asks for $1$ to be the identity element; the multiplication then extends by distributivity to all quaternions $q=x_0+x_1i+x_2j+x_3k \in \hh$ with $x_0,\ldots, x_3 \in \rr$. 

Since the beginning of last century there have been many attempts to determine a class of quaternion valued functions of one quaternionic variable playing the same role as the holomorphic functions of one complex variable. The most successful of such analogs is due to Fueter: in \cite{fueter2}, he defined a quaternionic function to be \emph{regular} if it solves the operator $\frac{\partial}{\partial \bar q} = \frac{1}{4} \left( \frac{\partial}{\partial x_0} + i \frac{\partial}{\partial x_1} + j \frac{\partial}{\partial x_2} +k \frac{\partial}{\partial x_3}\right).$
For an introduction to the theory of Fueter-regularity, see \cite{sudbery}. Over the years, this theory has been developed and generalized in many directions (see \cite{librodaniele,libroshapiro} and references therein). However, some of its features motivated the search for an alternative definition of regularity: for instance, the identity function is not Fueter-regular and the same holds for the most natural quaternionic polynomials. These are included in the class of \emph{quaternionic holomorphic} functions, defined by Fueter himself in \cite{fueter1} as solutions of the equation $\frac{\partial}{\partial \bar q} \Delta f(q) = 0,$
where $\Delta$ denotes the Laplacian in the four real variables $x_0,\ldots,x_3$. Notice, however,  that the class of quaternionic holomorphic functions is extremely large: it includes the class of harmonic functions of four real variables, which strictly includes that of Fueter-regular functions. 

An alternative definition of regularity was given in \cite{cullen} by Cullen, and developed by Gentili and Struppa in \cite{cras, advances}. It turned out that all polynomials and power series of the type
$$\sum_{n \in \nn} q^n a_n$$
with $a_n \in \hh$ define Cullen-regular functions on their sets of convergence, which are Euclidean balls centered at the origin of $\hh$. Conversely, every Cullen-regular function on such a ball admits a series expansion of this type. These two properties allow to prove that Cullen-regularity does not imply, nor is implied by Fueter-regularity: the function $q \mapsto q^2$ is an example of Cullen-regular function that is not harmonic, hence not Fueter-regular; the function $x_0 + i x_1 + j x_2 + k x_3 \mapsto x_0 + i x_1$ is Fueter-regular (hence harmonic and quaternionic holomorphic), but not Cullen-regular. After identifying $\hh = (\rr+i \rr)+(\rr+i \rr)j$ with $\cc^2$, the same examples prove that Cullen-regularity does not imply nor is implied by holomorphy in two complex variables. 

Quite recently, in \cite{advancesrevised}, Colombo, Gentili, Sabadini and Struppa identified a larger class of domains, called \emph{slice domains}, on which the study of regular functions is natural and not limited to quaternionic power series. This was truly a turning point in the theory, so that the term Cullen-regular function was substituted by \emph{slice regular function} (or simply \emph{regular function}) in this paper and in the subsequent works in this new context. Let us now review the definition and the basic properties of these functions. Let $\s = \{q \in \hh : q^2 = -1\}$ denote the 2-sphere of imaginary units of $\hh$, and recall that for all $I \in \s$ the subalgebra $\rr+I\rr$ is isomorphic to the complex field $\cc$.

\begin{definition}
Let $\Omega$ be a domain in $\hh$ and let $f : \Omega \to \hh$. For each $I \in \s$, let $\Omega_I = \Omega \cap (\rr + I \rr)$ and let $f_I = f_{|_{\Omega_I}}$ be the restriction of $f$ to the complex line $\rr + I \rr$. The restriction $f_I$ is called \emph{holomorphic} if it has continuous partial derivatives and
\begin{equation}
\bar \partial_I f(x+Iy) = \frac{1}{2} \left( \frac{\partial}{\partial x} + I \frac{\partial}{\partial y} \right) f_I(x+Iy) \equiv 0.
\end{equation}
The function $f$ is called \emph{(slice) regular} if, for all $I \in \s$, $f_I$ is holomorphic.
\end{definition}

The following Lemma clarifies the relationship between quaternionic regularity and complex holomorphy. 

\begin{lemma}[Splitting]\label{splitting}
Let $I \in \s$, let $\Omega_I$ be open in $\rr+I\rr$ and let $f_I : \Omega_I \to \hh$. The function $f_I$ is holomorphic if and only if, for all $J \in \s$ with $J \perp I$, there exist complex 
holomorphic functions 
$F,G : \Omega_I \to \rr+I\rr$ such that
\begin{equation}
f_I = F+GJ.
\end{equation} 
\end{lemma}

We already mentioned that any polynomial function of the type $q \mapsto a_0+qa_1+\ldots+q^n a_n$, with $a_l \in \hh$, is a regular function and that this class of examples extends as follows.

\begin{theorem}[Abel]\label{abel}
Let $\{a_n\}_n\in \nn$ be a sequence in $\hh$ and let $R \in [0,+\infty]$ be such that $1/R = \limsup_{n \in \nn} |a_n|^{1/n}$. The power series 
$
f(q) = \sum_{n \in \nn} q^n a_n
$ 
then converges absolutely and uniformly on compact sets in the Euclidean ball 
\begin{equation}
B(0,R) = \{q \in \hh : |q|<R\}, 
\end{equation}
where it defines a regular function $f$. On the other hand, $f(q)$ diverges at every point of $\hh \setminus \overline{B(0,R)}$.
\end{theorem}

For a proof, see \cite{advances}. The same paper showed that, conversely, all regular functions on $B(0,R)$ can be expressed as power series.

\begin{definition}
Let $f : \Omega \to \hh$ be a regular function. For each $I \in \s$, the \emph{$I$-derivative} of $f$ is defined as
$
\partial_I f(x+Iy) = \frac{1}{2} \left( \frac{\partial}{\partial x} - I \frac{\partial}{\partial y} \right) f_I(x+Iy)
$ 
on $\Omega_I$. The \emph{slice derivative} of $f$ is the regular function $\partial_s f : \Omega \to \hh$ defined to equal $\partial_I f$ on $\Omega_I$, for all $I \in \s$.
\end{definition}

\begin{theorem}[Series expansion]\label{seriesexpansion}
Let $R>0$ and let $f : B = B(0,R) \to \hh$ be a regular function. Then, for all $q \in B$,
$
f(q)=\sum_{n \in \nn} q^n \frac{1}{n!} f^{(n)}(0)
$ 
where $f^{(n)}$ denotes the $n$th slice derivative of $f$. In particular, $f \in C^{\infty}(B)$.
\end{theorem}

Theorems \ref{abel} and \ref{seriesexpansion} are fundamental in the study of regular quaternionic functions on balls $B=B(0,R)$ centered at the origin of $\hh$, since they identify these functions with power series. However, regular functions can have very different characteristics if we choose other types of domains of definition.

\begin{example}
Let $I \in \s$ and let $f: \hh \setminus \rr \to \hh$ be defined as follows:
$$f(q) = \left\{ 
\begin{array}{ll}
0 \ \mathrm{if} \ q \in \hh \setminus (\rr+I\rr)\\
1 \ \mathrm{if} \ q \in  (\rr+I\rr) \setminus \rr
\end{array}
\right.$$
This function is clearly regular, but not continuous.
\end{example}

The previous Example proves that if the domain $\Omega$ is not carefully chosen then a regular function $f : \Omega \to \hh$ does not even need to be continuous. It is possible to prevent such pathologies by imposing further conditions on the domain $\Omega$. The first such condition is described by next Definition (from \cite{advancesrevised}) and by the following Theorem (from \cite{poli}).

\begin{definition}
Let $\Omega$ be a domain in $\hh$ that intersects the real axis. $\Omega$ is called a \emph{slice domain} if, for all $I \in \s$, the intersection $\Omega_I$ with the complex line $\rr + I\rr$ is a domain of $\rr + I\rr$.
\end{definition}

\begin{theorem}[Identity Principle]\label{identity} 
Let $f,g$ be regular functions on a slice domain $\Omega$. If, for some $I \in \s$, $f$ and $g$ coincide on a subset of $\Omega_I$ having an accumulation point in $\Omega_I$, then $f = g$ in $\Omega$.
\end{theorem}

We now present a symmetry condition for the domains of definition, which guarantees continuity and differentiability for regular functions. In \cite{open} we noticed that the restriction of a power series $f(q) = \sum_{n \in \nn} q^n a_n$ to a sphere
$$x+y\s = \{x+Iy : I \in \s\}$$
is affine in the imaginary unit $I$, i.e. there exist $b,c \in \hh$ such that $f(x+Iy) = b+Ic$ for all $I \in \s$. As proven in \cite{cauchy,advancesrevised}, this is not only true for power series, but for all regular functions on the slice domains that have the following property.

\begin{definition}
A set $T \subseteq \hh$ is called \emph{axially symmetric} if, for all points $x+Iy \in T$ with $x,y \in \rr$ and $I \in \s$, the set $T$ contains the whole sphere $x+y\s$.
\end{definition}

Since no confusion can arise, we will refer to such a set as \emph{symmetric}, tout court. The following statement was proven in \cite{cauchy}, and it is a special case of a result proven in \cite{advancesrevised}.

\begin{theorem}[Representation Formula]
Let $f$ be a regular function on a symmetric slice domain $\Omega$ and let $x+y\s \subset \Omega$. For all $I,J \in \s$ 
\begin{eqnarray}
f(x+Jy) &=& \frac{1-JI}{2} f(x+Iy) +  \frac{1+JI}{2} f(x-Iy) \label{representationformula} \\ 
 &=& \frac{1}{2} \left[f(x+Iy) + f(x-Iy)\right] + \frac{JI}{2} \left[f(x-Iy) - f(x+Iy)\right].\nonumber
\end{eqnarray}
As a consequence, $f \in C^{\infty}(\Omega)$.
\end{theorem}

The result presented shows that a regular function on a symmetric slice domain $\Omega$ is uniquely determined by its restriction to a slice $\Omega_I$. This leads to the following Lemma, from \cite{cauchy} (a special case of the extension results proven in \cite{advancesrevised}).

\begin{lemma}\label{extensionlemma}
Let $\Omega$ be a symmetric slice domain and let $I \in \s$. If $f_I : \Omega_I \to \hh$ is holomorphic then there exists a unique regular function $g : \Omega \to \hh$ such that $g_I = f_I$ in $\Omega_I$. The function $g$ will be denoted by $\ext(f_I)$.
\end{lemma}

We conclude this first survey presenting quaternionic analogs of the domains of holomorphy.

\begin{definition}
The \emph{(axially) symmetric completion} of a set  $T \subseteq \hh$ is the smallest symmetric set $\widetilde{T}$ that contains $T$. In other words,
\begin{equation}
\widetilde{T} = \bigcup_{x+Iy \in T} (x+y\s).
\end{equation}
\end{definition}

\begin{theorem}[Extension]
Let $f$ be a regular function on a slice domain $\Omega$. There exists a unique regular function $\tilde f : \widetilde{\Omega} \to \hh$ that extends $f$ to the symmetric completion of its domain.
\end{theorem}

The present article studies the singularities of regular functions on symmetric slice domains. We conducted this study in \cite{poli} for the special case of Euclidean balls $B(0,R)$, but the approach there was somehow contrived because we were forced to only work with balls. The new context of symmetric slice domains allows a much more natural approach, as well as new results. We will point out the novelties along the exposition, but let us mention that we are now able to construct Laurent series at points other than the origin. Owing to the peculiarities of the non commutative setting, such a construction requires considerably different techniques than in the complex case. Let us also point out that we characterize essential singularities with a new result of Casorati-Weierstrass type.

Section \ref{sec:algebraic} presents the algebraic structure of regular functions: the set of regular functions on a symmetric slice domain $\Omega$ is a ring with respect to $+$ and to an appropriately defined multiplication, denoted by $*$. The multiplicative structure is (not trivially) related to the zero sets of these functions, presented in Section \ref{sec:zeros}. The results described in Sections \ref{sec:algebraic} and \ref{sec:zeros}, which already appeared in literature, are used extensively in our present study: as one may imagine, they are vital in the investigation of singularities.

The original part of the paper begins with Section \ref{sec:quotients}, which presents the algebraic properties of the quotients of regular functions on a symmetric slice domain $\Omega$. We indeed prove that the set of \emph{$*$-quotients} $f^{-*}*g$ of regular functions $f,g$ on $\Omega$ (with $f \not \equiv 0$) is the classical ring of quotients of the ring of regular functions on $\Omega$.

In Section \ref{sec:laurent} we study \emph{Laurent series}  
$$\sum_{n \in \zz}(q-p)^{*n} a_n$$ 
(where $(q-p)^{*n}$ denotes the $n$th power of $q \mapsto q-p$ with respect to $*$). The convergence of such a series is a delicate matter: in \cite{powerseries} we gave an estimate of $|(q-p)^{*n}|$ for $n \in \nn$ in terms of a non Euclidean distance $\sigma$ on $\hh$ (whose balls are pictured in Figure 1), an estimate which we presently improve to the optimal $|(q-p)^{*n}|\leq \sigma(q,p)^n$. Quite curiously, the corresponding estimate for the case $n<0$ is given in terms of a different function $\tau$. To conclude, we prove that the set of convergence of a Laurent series is of the type $\Sigma(p,R_1,R_2) = \{q \in \hh : \tau(q,p)>R_1, \sigma(q,p) <R_2\}$ and that the sum of the series is indeed a regular function in the interior of $\Sigma(p,R_1,R_2)$. Conversely, we prove that a regular function expands into Laurent series in each $\Sigma(p,R_1,R_2)$ contained in its domain of definition.

The expansion property allows us to classify the singularities of regular functions as \emph{removable}, \emph{essential} or \emph{poles} in Section \ref{sec:classification}. This classification is consistent with the one proposed in \cite{poli} in the special case of Euclidean balls $B(0,R)$. Furthermore, we define the quaternionic analogs of meromorphic functions, called \emph{semiregular} functions.

In Section \ref{sec:poles} we study poles in detail. We begin by proving that the quotient $f^{-*}*g$ of two regular functions on a symmetric slice domain $\Omega$ is semiregular in $\Omega$. We then prove the next result, where $ord_f(p)$ denotes the order of $p$ as a pole of $f$.

\begin{theorem}
Let $f$ be a semiregular function on a symmetric slice domain $\Omega$. Choose $p = x+yI \in \Omega$ and set $m = ord_f(p), n = ord_f(\bar p)$: without loss of generality $m \leq n$. There exist a neighborhood $U$ of $p$ in $\Omega$ that is a symmetric slice domain and a (unique) regular function $g : U \to \hh$ such that
\begin{eqnarray}
f(q) &=& [(q-p)^{*m}*(q-\bar p)^{*n}]^{-*}* g(q) = \\ 
 &=& \left[(q-x)^2+y^2\right]^{-n} (q-p)^{*(n-m)}* g(q) \nonumber
\end{eqnarray}
in $U \setminus (x+y\s)$. Moreover, if $n>0$ then $g(p) \neq 0, g(\bar p) \neq 0$.
\end{theorem}

We derive that the set of semiregular functions on a symmetric slice domain $\Omega$ is a division ring. This fact implies, in turn, that the function $g$ in the previous theorem extends to a semiregular function on $\Omega$. Finally, we study the distribution of the poles and their order: it turns out the relation between zeros and poles is more subtle than in the complex case.

Essential singularities are treated in Section \ref{sec:casorati}, where an analog of the Casorati-Weierstrass Theorem is proven
. The analogy to the complex case is complete when the essential singularity is located on the real axis, but not in the remaining cases, which are more elaborated.


\section{Algebraic prerequisites}\label{sec:algebraic}
In this Section we present the algebraic structure of regular functions. First of all, notice that the class of regular functions is endowed with an addition operation: if $f,g$ are regular functions on $\Omega$ then $f+g$ is regular in $\Omega$, too. The same does not hold for pointwise multiplication: $f \cdot g$ is not regular, except for some special cases. For this reason, in \cite{zeros} we used the following multiplicative operation: if $f(q) = \sum_{n \in \nn} q^n a_n, g(q) = \sum_{n \in \nn} q^n b_n$ are regular functions on $B(0,R)$ then setting
$$f*g(q) = \sum_{n \in \nn} q^n \sum_{k=0}^n a_k b_{n-k}$$
defines a regular function on the same ball $B(0,R)$. Notice that if $a_n \in \rr$ for all $n \in \nn$ then $f*g(q)=f(q)g(q)$. 
It turned out that the set of regular functions on a ball $B(0,R)$ is a ring with $+,*$. Recently, \cite{advancesrevised} extended this structure to all regular functions on symmetric slice domains. The construction of the $*$-product in this new setting relies upon the Splitting Lemma \ref{splitting} and upon Lemma \ref{extensionlemma}.

\begin{definition}
Let $f,g$ be regular functions on a symmetric slice domain $\Omega$. Choose $I,J \in \s$ with $I \perp J$ and let $F,G,H,K$ be holomorphic functions from $\Omega_I$ to $\rr+I\rr$ such that $f_I = F+GJ, g_I = H+KJ$. Consider the holomorphic function defined on $\Omega_I$ by
\begin{equation}
f_I*g_I(z) = \left[F(z)H(z)-G(z)\overline{K(\bar z)}\right] + \left[F(z)K(z)+G(z)\overline{H(\bar z)} \right]J.
\end{equation}
Its regular extension $\ext(f_I*g_I)$ to $\Omega$ is called the \emph{$*$-product} of $f$ and $g$ and it is denoted by $f*g$.
\end{definition}

It is possible to check directly that this Definition is coherent with the previous one in the special case $\Omega=B(0,R)$. Moreover, we can prove what follows.

\begin{proposition}
Let $\Omega$ be a symmetric slice domain. The definition of $*$-product is well posed  and the set of regular functions on $\Omega$ is a (non commutative) ring with respect to $+$ and $*$.
\end{proposition}

\begin{proof}
Let $f,g$ be regular functions on $\Omega$. By hypothesis, $\Omega$ intersects the real axis at some $r \in \rr$. The functions $q \mapsto f(q+r)$ and $q \mapsto g(q+r)$ are easily proven regular in $\Omega - r$, hence we may suppose without loss of generality $r=0$. This implies that there exists a ball $B=B(0,R) \subseteq \Omega$, with $R>0$. The restrictions $f_{|_B}$ and $g_{|_B}$ are power series.
We already observed that, for all $I \in \s$, $\ext(f_I*g_I)$ coincides with $f_{|_B}*g_{|_B}$ in $B$. In particular, for all $I,J \in \s$, $\ext(f_I*g_I)$ equals $\ext(f_J*g_J)$ in $B$. By the Identity Principle \ref{identity}, they coincide in $\Omega$. This proves that $f*g$ is well defined on $\Omega$.
The operation $*$ defined in this way is associative: $f*(g*h) = (f*g)*h$ because 
$$f_{|_B}*(g_{|_B} * h_{|_B}) = (f_{|_B}*g_{|_B})* h_{|_B}$$
The distributive law can be proven by the same technique. Finally, $*$ is clearly non commutative.
\end{proof}

As in the case of power series, the $*$-product coincides with the pointwise product for the following special class.

\begin{lemma}\label{realproduct}
Let $f,g$ be regular functions on a symmetric slice domain $\Omega$. If $f(\Omega_I)\subseteq \rr+I\rr$ for all $I \in \s$, then $fg$ is a regular function on $\Omega$ and $f*g= fg$.
\end{lemma}

In the special case $\Omega = B(0,R)$ for some $R>0$, Lemma \ref{realproduct} captures what we already observed for power series thanks to the fact that  $f(q) = \sum_{n \in \nn}q^na_n$ has real coefficients $a_n \in \rr$ if, and only if, $f(\Omega_I)\subseteq \rr+I\rr$ for all $I \in \s$. Moreover, there is an alternative expression of the $*$-product $f*g$. For power series, it is proven by direct computation (see \cite{zeros}), while the general case requires a different technique (see \cite{advancesrevised}). 

\begin{theorem}\label{formprod}
Let $f,g$ be regular functions on a symmetric slice domain $\Omega$.
For all $q\in\Omega$, if $f(q) = 0$ then $f*g(q) = 0$, else
\begin{equation}
f*g(q) = f(q)g(f(q)^{-1} q f(q)).
\end{equation}
\end{theorem}

Let us now present two other operations introduced in \cite{zeros} for a regular function $f(q) = \sum_{n \in \nn} q^n a_n$ on a ball $B(0,R)$: the \emph{regular conjugate} of $f$ is the regular function defined by
$$f^c(q) = \sum_{n \in \nn} q^n \bar a_n$$
on the same ball $B(0,R)$, the \emph{symmetrization} of $f$ is the function
$$f^s = f*f^c = f^c*f.$$
These operations were defined in order to study the zero set, but they also allowed us to construct the ring of quotients of regular functions on $B(0,R)$ in \cite{poli}. In \cite{advancesrevised}, they were extended to all symmetric slice domains in the following manner.

\begin{definition}
Let $f$ be a regular function on a symmetric slice domain $\Omega$. Choose $I,J \in \s$ with $I \perp J$ and let $F,G$ be holomorphic functions from $\Omega_I$ to $\rr+I\rr$ such that $f_I = F+GJ$. If $f_I^c$ is the holomorphic function defined on $\Omega_I$ by
\begin{equation}
f_I^c(z) = \overline{F(\bar z)} - G(z)J.
\end{equation}
then the \emph{regular conjugate} of $f$ is the regular function defined on $\Omega$ by $f^c=\ext(f_I^c)$. The \emph{symmetrization} of $f$ is the regular function defined on $\Omega$ by $f^s = f*f^c = f^c*f$.
\end{definition}

This Definition is coherent with the previous one in the special case $\Omega=B(0,R)$ (by direct computation). It is also possible to prove that the Definition is well posed using the same technique as in the case of the $*$-product. Notice, moreover, that if $f(\Omega_I)\subseteq \rr+I\rr$ for all $I \in \s$, then $f^c(q)=f(q)$ and $f^s(q) = f(q)^2$ for all $q \in \Omega$.


\section{Prerequisites on the zero sets}\label{sec:zeros}

As one may imagine, the study of singularities is intimately related to that of zeros. In order to make our presentation self-contained, we now survey some properties of the zero sets of regular functions which are known in literature.  

The study of zero sets began with the proof, in \cite{advances}, of the following result (in the special case $\Omega=B(0,R)$).

\begin{theorem}\label{affinezeros}
If $f$ is a regular function on a symmetric slice domain $\Omega$ vanishing at a point $x+Jy$ then either $f$ vanishes identically in $x+y\s$ or $f$ does not have any other zero in $x+y\s$. The first case always applies if $f(\Omega_I) \subseteq \rr+I\rr$ for all $I \in \s$.
\end{theorem}

Theorem \ref{affinezeros} was proven in its most general statement in \cite{advancesrevised}, as a consequence of the fact that $f_{|_{x+y\s}}$ is affine. The study of the zero sets continued with their algebraic properties, proven in \cite{zeros, milan} for power series and polynomials, and in \cite{advancesrevised,zerosopen} for regular functions on symmetric slice domains. First consider this consequence of Theorem \ref{formprod}.

\begin{corollary}
If $f,g$ are regular functions on a symmetric slice domain $\Omega$ and $q \in \Omega$, then $f*g(q) = 0$ if and only if either $f(q) = 0$ or $g(f(q)^{-1} q f(q))=0$. As a consequence, for every zero of $f*g$ in a sphere $x+y\s$ there exists a zero of $f$ or a zero of $g$ in $x+y\s$.
\end{corollary}

However, if we denote the zero set of $f$ as $Z_f$, there need not be a one-to-one correspondence between $Z_{f*g}$ and $Z_{f}\cup Z_{g}$, as proven by the following polynomial case.

\begin{proposition}\label{roots}
Let $\alpha,\beta\in \hh$ and $P(q)=(q-\alpha)*(q-\beta)$. 
\begin{enumerate}
\item If $\beta$ does not lie in the same sphere $x+y\s$ as $\alpha$ then $P$ has two zeros, $\alpha$ and $(\beta - \bar \alpha)^{-1}\beta(\beta - \bar \alpha)$. 
\item If $\alpha,\beta$ lie in the same sphere $x+y\s$ but $\alpha \neq \bar \beta$ then $P$ only vanishes at $\alpha$. 
\item Finally, if $\alpha=\bar \beta \in x+y\s$ then the zero set of $P$ is $x+y\s$.
\end{enumerate}
\end{proposition}

It is also possible to study the effect of conjugation and symmetrization on the zero set.

\begin{lemma}\label{symmetrizationreal}
Let $f$ be a regular function on a symmetric slice domain $\Omega$ and let $f^s$ be its symmetrization. Then $f^s(\Omega_I) \subseteq \rr+I\rr$ for all $I \in \s$. In particular, for each $S=x+y\s \subset \Omega$ either $f^s$ vanishes identically in $S$ or it has no zeros in $S$.
\end{lemma}

\begin{proposition}\label{conjugatezeros}
Let $f$ be a regular function on a symmetric slice domain $\Omega$ and choose $S=x+y\s\subset \Omega$. The zeros of $f$ in $S$ are in one-to-one correspondence with those of $f^c$. Furthermore, $f^s$ vanishes identically on $S$ if and only if $f^s$ has a zero in $S$, if and only if $f$ has a zero in $S$ (if and only if $f^c$ has a zero in $S$).
\end{proposition}


The algebraic properties of the zeros presented so far allowed the following topological characterization.

\begin{theorem}[Structure of the Zero Set]
Let $f$ be a regular function on a symmetric slice domain $\Omega$. If $f$ does not vanish identically, then the zero set of $f$ consists of isolated points or isolated 2-spheres of the form $x+y\mathbb{S}$.
\end{theorem}


Before concluding this Section, let us speak about factorization and multiplicity. Each zero of a regular function can be ``factored out'': we proved it for power series in \cite{zeros} and we easily extend the result as follows.

\begin{proposition}
Let $f$ be a regular function on a symmetric slice domain $\Omega$. A point $p \in \Omega$ is a zero of $f$ if and only if there exists a regular function $g : \Omega \to \hh$ such that
$
f(q) = (q-p)*g(q).
$
\end{proposition}

This naturally leads to the following concept of multiplicity (defined in \cite{zeros} for power series).

\begin{definition}
Let $f$ be a regular function on a symmetric slice domain $\Omega$ and let $p \in \Omega$. We define the \emph{(classical) multiplicity} of $p$ as a zero of $f$ and denote by $m_f(p)$ the largest $n \in \nn$ such that
$
f(q)=(q-p)^{*n}*g(q)
$ 
for some regular $g: \Omega \to \hh$.
\end{definition}

Notice that if $p \in \rr+I\rr$, then $f(q) = (q-p)^{*n}*g(q)$ if and only if $f_I(z) = (z-p)^n g_I(z)$. Hence the classical multiplicity of $p$ as a zero of $f$ coincides with the multiplicity of $p$ as a zero of the holomorphic function $f_I$, intended as the largest $n \in \nn$ such that $f_I(z) = (z-p)^n g_I(z)$. Notice that, if $f_I$ splits as $f_I=F+GJ$ according to Lemma \ref{splitting}, then the multiplicity of $f_I$ at $p$ is the minimum between the multiplicity of $F$ and that of $G$ at $p$. This proves a posteriori that the definition is well posed.

Even though the classical multiplicity is a consistent generalization of complex multiplicity, it does not lead to analogous results for polynomials. Indeed, with this Definition, a polynomial of finite degree can have infinitely many roots with positive multiplicity.

\begin{example}\label{sphere}
The polynomial $P(q) = q^2+1 = (q-I)*(q+I)$ has multiplicity $m_P(I) = 1$ at all $I \in \s$.
\end{example}

Even if we consider a polynomial having only isolated (hence finitely many) roots, we still cannot relate their number to the degree of the polynomial. Indeed, the degree of a polynomial can exceed the sum of the multiplicities of its zeros.

\begin{example}\label{oneroot}
Take $I,J \in \s$ with $I \neq \pm J$ and let
$$P(q) = (q-I)*(q-J) = q^2-q(I+J)+IJ.$$
By Proposition \ref{roots}, the zero set of $P$ is $\{I\}$. Notice that $m_P(I)=1$, while $f$ has degree $2$.
\end{example}

For this reason \cite{milan} introduced alternative notions of multiplicity for the roots of polynomials. Let us present them by means of the following factorization result. For the purpose of this paper, it is useful to extend both the result and the Definition from polynomials to regular functions on symmetric slice domains.

\begin{theorem}\label{factorization}
Let $f$ be a regular function on a symmetric slice domain $\Omega$, suppose $f \not \equiv 0$ and let $x+y\s \subset\Omega$. There exist $m \in \nn$ and a regular function $\tilde f:\Omega \to \hh$ not identically zero in $x+y\s$ such that
\begin{equation}
f(q) = [(q-x)^2+y^2]^m \tilde f(q).
\end{equation}
If $\tilde f$ has a zero $p_1 \in x+y\s$ then such a zero is unique and there exist $n \in \nn$, $p_2,...,p_n \in x+y\s$ (with $p_i \neq \bar p_{i+1}$ for all $i \in \{1,\ldots,n-1\}$) such that
\begin{equation}
\tilde f(q) = (q-p_1)*(q-p_2)*...*(q-p_n)*g(q)
\end{equation}
for some regular function $g:\Omega \to \hh$ which does not have zeros in $x+y\s$.
\end{theorem}

\begin{definition}
In the situation of Theorem \ref{factorization}, we say that $f$ has \emph{spherical multiplicity} $2m$ at $x+y\s$ and that $f$ has \emph{isolated multiplicity} $n$ at $p_1$.
\end{definition}

As observed in \cite{milan}, the degree of a polynomial equals the sum of the spherical multiplicities and of the isolated multiplicities of its zeros. For instance: in the case of Example \ref{sphere}, $P$ has spherical multiplicity $2$ at $\s$; in Example \ref{oneroot}, $P$ has spherical multiplicity $0$ at $\s$ and isolated multiplicity $2$ at $I$.



\section{Quotients of regular functions}\label{sec:quotients}

Our present study regards the singularities of regular functions and (as in the complex case) it requires the construction of the ring of quotients of these functions. Let us begin by presenting the reciprocal of $f$ with respect to $*$, defined in \cite{open} for power series and extended in \cite{advancesrevised} as follows.

\begin{definition}
Let $f$ be a regular function on a symmetric slice domain $\Omega$. If $f\not \equiv0$, the \emph{regular reciprocal} of $f$ is the function defined on $\Omega \setminus Z_{f^s}$ by
\begin{equation}
f^{-*}=\frac{1}{f^s}f^c.
\end{equation}
\end{definition}

In order to prove the regularity of $f^{-*}$, we will make use of the next Lemma.

\begin{lemma}
Let $h$ be a regular function on a symmetric slice domain $\Omega$. If $h\not \equiv0$ and if $h(\Omega_I)\subseteq \rr+I\rr$ for all $I \in \s$ then $\frac{1}{h}$ is a regular function in $\widehat{\Omega} = \Omega \setminus Z_h$, mapping $\widehat{\Omega}_I$ to $\rr+I\rr$ for all $I \in \s$.
\end{lemma}

\begin{proof}
Set $g(q)= \frac{1}{h(q)}$ for all $q \in \widehat{\Omega}=\Omega \setminus Z_h$. We will prove that $g$ is regular by observing that, for all $I \in \s$, the restriction $g_I$ is holomorphic. Let $I \in \s$. Chosen $J \in \s$ with $J \perp I$, the Splitting Lemma \ref{splitting} guarantees the existence of holomorphic functions $F,G : \Omega_I \to \rr+I\rr$ such that $h_I = F+GJ$. Since $h(\Omega_I)\subseteq \rr+I\rr$, we must have $G \equiv 0$ and $h_I = F$. Hence 
$$g_I(z) = \frac{1}{h_I(z)} = \frac{1}{F(z)}$$ 
for all $z \in \widehat{\Omega}_I$ and $g_I$ is a holomorphic function from $\widehat{\Omega}_I$ to $\rr+I\rr$.
\end{proof}

\begin{theorem}
Let $f$ be a regular function on a symmetric slice domain $\Omega$. The function $f^{-*}$ is regular in $\Omega \setminus Z_{f^s}$, which is a symmetric slice domain, and
$
f*f^{-*} = f^{-*}*f = 1
$ 
in $\Omega \setminus Z_{f^s}$.
\end{theorem}

\begin{proof}
The regularity of $f^{-*} = \frac{1}{f^s}f^c$ is proven as follows: the previous Lemma applies to $h=f^s$ thanks to Lemma \ref{symmetrizationreal}; the regularity of $\frac{1}{f^s}f^c$ then follows by Lemma \ref{realproduct}. Furthermore, $\Omega \setminus Z_{f^s}$ is still a symmetric slice domain since $f \not \equiv 0$ implies that $Z_{f^s}$ consists of isolated real points and isolated 2-spheres $x+y\s$. Finally,
$$f^{-*} *f = \frac{1}{f^s} f^c*f =  \frac{1}{f^s} f^s =1$$
and, similarly, $f*f^{-*} = 1$.
\end{proof} 

Let us now present the ring of quotients of regular functions.
Recall that, in the complex case, the set of quotients $\frac{F}{G}$ of holomorphic functions $F,G$ on a disk $\Delta$ becomes a field when endowed with the usual operations of addition and multiplication. More precisely, it is the field of quotients of the integral domain (the commutative ring with no zero divisors) obtained by endowing the set of holomorphic functions $F$ on $\Delta$ with the natural addition and multiplication. As explained in \cite{rowen} (see also \cite{cohn, lam2}), the concept of field of quotients of an integral domain can be generalized to the non-commutative case. 

\begin{theorem}
Define a \emph{left Ore domain} as a domain (a ring with no zero divisors) $(D,+,\cdot)$ such that $Da \cap Db \neq \{0\}$ for all $a,b \in D \setminus \{0\}$. If this is the case, then the set of formal quotients $L = \{a^{-1}b : a,b \in D\}$ can be endowed with operations $+, \cdot$ such that: 
\begin{enumerate}
\item[(1)] $D$ is isomorphic to a subring  of $L$ (namely $ \{1^{-1}a : a \in D\}$); 
\item[(2)] $L$ is a skew field, i.e. a ring where every non-zero element has a multiplicative inverse (namely, $(a^{-1}b)^{-1} = b^{-1}a$).
\end{enumerate}
The ring $L$ is called the \emph{classical left ring of quotients} of $D$ and, up to isomorphism, it is the only ring having the properties (1) and (2). 
\end{theorem}
On a \emph{right Ore domain} $D$, defined by $aD \cap bD \neq \{0\}$ for all $a,b \in D \setminus \{0\}$, the \emph{classical right ring of quotients} is similarly constructed. If $D$ is both a left and a right Ore domain, then (by the uniqueness property) the two rings of quotients are isomorphic and we speak of the \emph{classical ring of quotients} of $D$. 
We now construct the ring of quotients of regular functions on a symmetric slice domain.

\begin{theorem}\label{quotients}
Let $\Omega$ be a symmetric slice domain. The set of \emph{$*$-quotients}
\begin{equation}
\mathcal{L}(\Omega) = \{f^{-*}*g : f,g \mathrm{\ regular\ in\ } \Omega, f \not \equiv 0\}
\end{equation}
is a division ring with respect to $+, *$. Moreover, the ring of regular functions on $\Omega$ is a left and right Ore domain and $\mathcal{L}(\Omega)$ is its classical ring of quotients.
\end{theorem}

\begin{proof}
For all regular $f,g,h,k : \Omega \to \hh$ with $f,h \not \equiv 0$, the sum and $*$-product of $f^{-*}*g$ and $h^{-*}*g$ are regular functions on $\Omega \setminus (Z_{f^s}\cup Z_{h^s})$ that prove to be elements of $\mathcal{L}(\Omega)$ by means of the following formulae:
$$(f^{-*}*g)+(h^{-*}*k) = \frac{1}{f^sh^s} (h^s f^c*g+ f^s h^c*k),$$
$$(f^{-*}*g)*(h^{-*}*k) = \frac{1}{f^sh^s} f^c*g*h^c*k.$$
We easily derive that $\mathcal{L}(\Omega)$ is a ring with respect to $+,*$. Furthermore, $\mathcal{L}(\Omega)$ is a division ring since
$$(f^{-*}*g)*(g^{-*}*f) = \frac{1}{f^sg^s} f^c*g*g^c*f = \frac{1}{f^sg^s} f^c*g^s*f = \frac{1}{f^sg^s} f^s g^s = 1.$$

The ring $D$ of regular functions on $\Omega$ is a domain, since $f*g \equiv 0$ if and only if $f \equiv 0$ or $g\equiv 0$. Moreover, $D$ is a left Ore domain: if $f,g \not \equiv 0$ then $(D*f) \cap (D*g)$ contains the non-zero element $f^s g^s = g^s f^s$, which can be obtained as $(g^s * f^c )*f$ or as $(f^s*g^c)*g$. Similarly,  $D$ is a right Ore domain. Thus the classical ring of quotients of $D$ is well defined. It must be isomorphic to $\mathcal{L}(\Omega)$ by the uniqueness property: $\mathcal{L}(\Omega)$ is a division ring having $D$ as a subring and the inclusion $D \to \mathcal{L}(\Omega) \ \ f \mapsto f = 1^{-*} *f$ is a ring homomorphism.
\end{proof}

We had already proven Theorem \ref{quotients} for the special case $\Omega = B(0,R)$ in \cite{poli}, but our approach has changed: in that paper, starting with multiplication $*$ only defined on power series, we had to define $*$ on quotients in an abstract manner and prove its consistency; we now start with $(f^{-*}*g)*(h^{-*}*g)$ already defined (since $\Omega \setminus (Z_{f^s}\cup Z_{h^s})$ is a symmetric slice domain when $\Omega$ is) and prove that this $*$-product is itself a quotient.
This construction is an important tool in the study of the singularities of regular functions, which we will soon present. Before doing so, we introduce another tool: Laurent series are studied in the next Section.


\section{Laurent series and expansions}\label{sec:laurent}

We will now study Laurent series and expansions. In \cite{poli} we proved that a quaternionic Laurent series centered at $0$
$$\sum_{n \in \zz} q^n a_n = \sum_{n\geq0} q^n a_n + \sum_{m>0} q^{-m} a_{-m}$$
converges absolutely and uniformly on compact sets in a spherical shell $A = \{q \in \hh : R_1< |q| < R_2\},$
where its sum defines a regular function. The shell can be chosen so that the series diverges outside $\overline{A}$. In other words, the analog of Theorem \ref{abel} for Laurent series holds. Conversely, an appropriate analog of Theorem \ref{seriesexpansion} was proven in the same paper. These results, however, are not as significant as one may imagine: they still hold if we move the center of the series from $0$ to another point of the real axis, but not if we center it a point $p \in \hh \setminus \rr$. 

\begin{example}
The function $q \mapsto (q-p)^2=q^2-qp-pq+p^2$ is not regular, unless $p \in \rr$ (in which case $pq=qp$). The same holds for $(q-p)^n$ for all $n \in \zz \setminus \{0,1\}$.
\end{example}

For this reason, we are encouraged to give the following Definition. For all $p \in \hh$ and $n \in \nn$, we denote by $(q-p)^{*n} = (q-p)*\ldots*(q-p)$ the $n$th $*$-power of $q \mapsto q-p$; for all $m \in \nn$, let us write $(q-p)^{*(-m)} = (q-p)^{-*m}$ for the regular reciprocal of $(q-p)^{*m}$.

\begin{definition}
Let $p \in \hh$. For any $\{a_n\}_{n\in \zz}$ in $\hh$ we call
$
\sum_{n \in \zz}(q-p)^{*n} a_n
$ 
the \emph{Laurent series centered at $p$} associated to $\{a_n\}_{n\in \zz}$.
\end{definition}

Clearly, every addend in the series 
is a regular function. Let us now study its convergence. In \cite{powerseries} we studied \emph{regular power series} $\sum_{n \in \nn}(q-p)^{*n} a_n.$ In order to study their convergence, we estimated $|(q-p)^{*n}|$ as $2 \sigma(q,p)^n$ for all $n \in \nn$, where $\sigma$ is defined as follows.

\begin{definition}
For all $p,q \in \hh$
\begin{equation}
\sigma(q,p) = \left\{
\begin{array}{ll}
|q-p| & \mathrm{if\ } p,q \mathrm{\ lie\ on\ the\ same\ complex\ line\ } \rr+I\rr\\
\omega(q,p) &  \mathrm{otherwise}
\end{array}
\right.
\end{equation}
where
\begin{equation}
\omega(q,p) = \sqrt{\left[Re(q)-Re(p)\right]^2 + \left[|Im(q)| + |Im(p)|\right]^2}.
\end{equation}
\end{definition}

In other words, if $p \in \rr+I\rr$ then for $q \in \rr+I\rr$ $\sigma(q,p) = |q-p|$, while for $q \in \hh \setminus (\rr+I\rr)$ we have 
$$\sigma(q,p) = \max\{|z-p|, |\bar z-p|\}$$
where $z,\bar z$ are the points of $\rr+I\rr$ such that $Re(z) =Re(\bar z) = Re(q)$ and $|Im(z)| =|Im(\bar z)|= |Im(q)|$.
In turned out that

\begin{proposition}
The function $\sigma : \hh \times \hh \to \rr$ is a distance.
\end{proposition}

However, $\sigma$ is not topologically equivalent to the Euclidean distance, as proven by the following discussion.

\begin{remark} Let $p \in \rr+I\rr \subset \hh$ and $R \geq 0$. Denote the $\sigma$-ball of radius $R$ centered at $p$ as
\begin{equation}
\Sigma(p,R) = \{q \in \hh : \sigma(q,p)<R\}.
\end{equation}
Furthermore, define $\Omega(p,R) =\{q \in \hh : \omega(q,p) < R\}$.
\begin{enumerate}
\item If $R \leq |Im(p)|$ then $\Sigma(p,R)$ reduces to the Euclidean disc $\Delta_I(p,R)$ in $\rr+I\rr$ (and the set $\Omega(p,R)$ is empty).
\item If $R > |Im(p)| > 0$ then $\Sigma(p,R)= \Delta_I(p,R) \cup \Omega(p,R)$ and $\Omega(p,R)$ is the symmetric completion of $\Delta_I(p,R) \cap \Delta_I(\bar p,R)$.
\item If $p \in \rr$, i.e. $Im(p)=0$ then $\Sigma(p,R)$ coincides with the Euclidean ball $B(p,R)$ (and with $\Omega(p,R)$).
\end{enumerate}
\end{remark}
\begin{figure}[h]
  \begin{center}
  \fbox{\includegraphics[height=7cm]{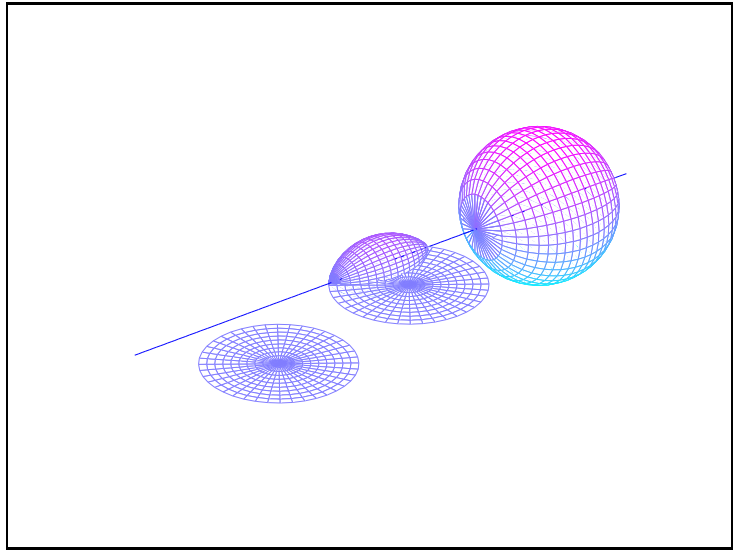}}
\end{center}
  \caption{A view in $\rr+i\rr+j\rr$ of $\sigma$-balls $\Sigma(p,R)$ centered at points $p \in \rr+i\rr$ and having $|Im(p)| \geq R$, $0 < |Im(p)| <R$ and $Im(p) = 0$, respectively.}
\end{figure}
In all three cases, the interior of $\Sigma(p,R)$ with respect to the Euclidean topology is $\Omega(p,R)$ and its Euclidean closure is $\overline{\Sigma(p,R)} =  \{q \in \hh : \sigma(q,p) \leq R\}$. Notice that $\Omega(p,R)$ is a symmetric slice domain when it is not empty. Furthermore, $p$ is in the interior of $\Sigma(p,R)$, i.e. $p \in \Omega(p,R)$, if and only if $2 |Im(p)| < R$.

By means of the estimate $|(q-p)^{*n}| \leq 2 \sigma(q,p)^n$ and of the root test, we proved the ``$\sigma$-analog'' of the Abel Theorem \ref{abel}.

\begin{theorem}
Choose any sequence $\{a_n\}_{n \in \nn}$ in $\hh$ and let $R \in [0, +\infty]$ be such that $1/R = \limsup_{n \to +\infty} |a_n|^{1/n}$. For any fixed $p \in \hh$, the series 
$
f(q) = \sum_{n \in \nn} (q-p)^{*n} a_n
$ 
converges absolutely and uniformly on the compact subsets of $\Sigma(p,R)$ and it does not converge at any point of $\hh \setminus \overline{\Sigma(p,R)}$ (we call $R$ the \emph{$\sigma$-radius of convergence} of $f(q)$). Furthermore, if $\Omega(p,R) \neq \emptyset$ then the sum of the series defines a regular function $f : \Omega(p,R) \to \hh$.
\end{theorem}

Conversely, we proved that a regular function on a symmetric slice domain $\Omega$ expands into regular power series at each point $p \in \Omega$. We thus defined and studied appropriate notions of analyticity (the interested reader is referred to \cite{powerseries}). 

We now study the convergence of Laurent series, beginning with next Definition.

\begin{definition}
For all $p,q \in \hh$ we define
\begin{equation}
\tau(q,p) = \left\{
\begin{array}{l}
|q-p|  \mathrm{\ if\ } p,q \mathrm{\ lie\ on\ the\ same\ complex\ line\ } \rr+I\rr\\
\sqrt{\left[Re(q)-Re(p)\right]^2 + \left[|Im(q)| - |Im(p)|\right]^2}   \mathrm{\ otherwise}
\end{array}
\right.
\end{equation}
Equivalently: if $p \in \rr+I\rr$ then for all $q \in \rr+I\rr$ we set $\tau(q,p) = |q-p|$; for all $q \in \hh \setminus (\rr+I\rr)$, we instead set $\tau(q,p) = \min\{|z-p|, |\bar z-p|\}$ where $z, \bar z$ are the points of $\rr+I\rr$ having the same real part and the same modulus of the imaginary part as $q$.
\end{definition}

We improve the estimate of $|(q-p)^{*n}|$ given in \cite{powerseries}, and find the corresponding estimate for the case $n<0$.

\begin{theorem}
Fix $p \in \hh$. Then $|(q-p)^{*n}| \leq \sigma(q,p)^n$ for all $n \in \nn$ and $|(q-p)^{-*m}| \leq \tau(q,p)^{-m}$ for all $m \in \nn$. Moreover,
\begin{eqnarray}
\lim_{n \to +\infty} |(q-p)^{*n}|^{1/n} &=& \sigma(q,p)\\
\lim_{m \to +\infty} |(q-p)^{-*m}|^{1/m}&=& \frac{1}{\tau(q,p)}.\nonumber
\end{eqnarray}
\end{theorem}

\begin{proof}
Let $\rr+I\rr$ be the complex line through $p$. For all $z = x+Iy \in \rr+I\rr$, $z$ and $p$ commute so that $(q-p)^{*n}$ equals $(z-p)^n$ when computed at $q=z$. For any $q = x+Jy$, Formula (\ref{representationformula}) implies:
$$(q-p)^{*n} = \frac{1}{2} \left[(z-p)^n+(\bar z-p)^n\right] - \frac{JI}{2} \left[(z-p)^n-(\bar z-p)^n\right].$$
Notice that, in vector notation,  $JI = - <J,I> + J\times I$. If $<J,I> = \cos\theta$ then $J\times I = \sin\theta\ L$ where $L \in \s$ is orthogonal to $I$. Hence
$$|(q-p)^{*n}|^2 = \left|\frac{1}{2}[(z-p)^n + (\bar z-p)^n] + \frac{\cos\theta}{2} [(z-p)^n -(\bar z-p)^n]\right|^2 +$$
$$+\frac{\sin^2\theta}{4} \left|(z-p)^n - (\bar z-p)^n\right|^2=$$
$$=\frac{1}{4}\left|(z-p)^n + (\bar z-p)^n\right|^2+\frac{\cos^2\theta + \sin^2\theta}{4} \left|(z-p)^n - (\bar z-p)^n\right|^2+$$
$$ + \frac{\cos\theta}{2}\left<(z-p)^n + (\bar z-p)^n, (z-p)^n -(\bar z-p)^n\right>$$
which attains its maximum value when $\cos \theta = 1$ or when $\cos \theta = -1$, in other words when $J = I$ or when $J = -I$, i.e. at $q = z$ or at $q=\bar z$. By analogous computations, $|(q-p)^{-*m}|$ attains its maximum value at $q = z$ or at $q=\bar z$.
Hence, for all $q \not \in \rr+I\rr$,
$$|(q-p)^{*n}| \leq \max\{|z-p|^n, |\bar z-p|^n\} = \max\{|z-p|, |\bar z-p|\}^n = \sigma(q,p)^n,$$
$$|(q-p)^{-*m}| \leq \max \left\{\frac{1}{|z-p|^{m}}, \frac{1}{|\bar z-p|^{m}} \right\} =$$
$$= \frac{1}{\min \left\{|z-p|, |\bar z-p|\right\}^{m}} =\frac{1}{\tau(q,p)^{m}}.$$
If, on the contrary, $q \in \rr+I\rr$ then $|(q-p)^{*n}| = |q-p|^n = \sigma(q,p)^n$ and $|(q-p)^{-*m}| = \frac{1}{|q-p|^{m}} = \frac{1}{\tau(q,p)^{m}}$. This proves the first statement.

Let us now prove the second statement. It is trivial in the case $q \in \rr+I\rr$. When $q \not \in \rr+I\rr$, i.e. $y \neq 0$ and $J \neq \pm I$, we prove it as follows. If $|z-p|=|\bar z-p|$ then $(y-|Im(p)|)^2 = (-y - |Im(p)|)^2$, hence $Im(p) = 0$, i.e. $p \in \rr$, so that $(q-p)^{*n}=(q-p)^n$, $(q-p)^{-*m}=(q-p)^{-m}$ and the thesis is trivial. Let us thus suppose $|z-p|\neq |\bar z-p|$. Without loss of generality, $|z-p|<|\bar z-p|$ and in particular $\sigma(q,p)= \max\{|z-p|,|\bar z-p|\} = |\bar z-p|$, $\tau(q,p)= \min\{|z-p|, |\bar z-p|\} = |z-p|$. Since
$$(q-p)^{*n} = \frac{1-JI}{2} (z-p)^n+ \frac{1+JI}{2} (\bar z-p)^n,$$
$$(q-p)^{-*m} = \frac{1-JI}{2} \frac{1}{(z-p)^{m}}+ \frac{1+JI}{2} \frac{1}{(\bar z-p)^{m}}$$
we have
$$\frac{\left|(q-p)^{*n}\right|}{\sigma(q,p)^n} = \left|\frac{1-JI}{2} \left(\frac{z-p}{\bar z-p}\right)^n+ \frac{1+JI}{2}\right|,$$
$${\left|(q-p)^{-*m}\right|}{\tau(q,p)^{m}} = \left|\frac{1-JI}{2} + \frac{1+JI}{2} \left(\frac{z-p}{\bar z-p}\right)^{m}\right|.$$
where $\left|\frac{z-p}{\bar z-p}\right|<1$ and $\frac{1\pm JI}{2}\neq 0$. We easily conclude that 
$$\lim_{n \to +\infty} \frac{|(q-p)^{*n}|^{1/n}} {\sigma(q,p)} =1,$$
$$\lim_{m \to +\infty} |(q-p)^{-*m}|^{1/m} \tau(q,p)=1$$
which is equivalent to our thesis.
\end{proof}

The function $\tau$ is clearly not a distance, since $\tau(q,p) = 0$ for all $q$ lying in the same 2-sphere $x+y\s$ as $p$ except $\bar p$. However, it is useful for our purposes to study its level sets.

\begin{remark} Chosen $p \in \hh$ and $R \geq 0$, let us denote $T(p,R) = \{q \in \hh : \tau(q,p) < R\}$. If $p = x+Iy$ then 
$$T(p,R) = \bigcup_{J\neq -I} \Delta_J(x+Jy,R).$$
\begin{enumerate}
\item If $R \leq |Im(p)|$ then $T(p,R)$ is the symmetric completion of the Euclidean disc $\Delta_I(p,R)\subset \rr+I\rr$, minus the Euclidean disc $\Delta_I(\bar p,R)\subset \rr+I\rr$.
\item If $R > |Im(p)| > 0$ then $T(p,R)$ is the symmetric completion of the Euclidean disc $\Delta_I(p,R)\subset \rr+I\rr$, minus $\Delta_I(\bar p,R) \setminus \Delta_I(p,R)$.
\item If $p \in \rr$, i.e. $Im(p)=0$ then $T(p,R)$ coincides with the Euclidean ball $B(p,R)$.
\end{enumerate}
\end{remark}

In all cases $T(p,R)$ is open in the Euclidean topology. Notice that 
$$\{q \in \hh : \tau(q,p) \leq R\} = \bigcup_{J\neq -I} \overline{\Delta_J(x+Jy,R)}$$ 
is not, in general, the closure of $T(p,R)$ with respect to the Euclidean topology. We instead have:
$$\overline{T(p,R)} =  \bigcup_{J\in \s} \overline{\Delta_J(x+Jy,R)}.$$

\begin{remark}
With respect to the Euclidean topology, $\sigma : \hh \times \hh \to \rr$ is lower semicontinuous and $\tau : \hh \times \hh \to \rr$ is upper semicontinuous.
\end{remark}

We now define the sets which will prove to be the sets of convergence of Laurent series.

\begin{definition}
For $0\leq R_1 <R_2\leq +\infty$, we define
\begin{equation}
\Sigma(p,R_1,R_2) = \{q \in \hh : \tau(q,p)>R_1, \sigma(q,p) <R_2\},
\end{equation}
\begin{equation}
\Omega(p,R_1,R_2) = \Omega(p,R_2)\setminus\overline{T(p,R_1)},
\end{equation}
where $\Omega(p,R_2) = \{q \in \hh : \omega(q,p) <R_2\}$.
\end{definition}

Notice that the interior of $\Sigma(p,R_1,R_2)$ is $\Omega(p,R_1,R_2)$. On the other hand,
$$\overline{\Sigma(p,R_1,R_2)} = \{q \in \hh : \tau(q,p)\geq R_1, \sigma(q,p) \leq R_2\}$$
as a consequence of the previous Remark. We are now ready to study the convergence of Laurent series.

\begin{theorem}
Choose any sequence $\{a_n\}_{n \in \zz}$ in $\hh$. Let $R_1,R_2 \in [0, +\infty]$ be such that $R_1 = \limsup_{m \to +\infty} |a_{-m}|^{1/m}, 1/R_2 = \limsup_{n \to +\infty} |a_n|^{1/n}$.
For all $p \in \hh$ the series 
\begin{equation}
f(q) = \sum_{n \in \zz} (q-p)^{*n} a_n
\end{equation}
converges absolutely and uniformly on the compact subsets of $\Sigma(p,R_1,R_2)$ and it does not converge at any point of $T(p,R_1)$ nor at any point of $\hh \setminus \overline{\Sigma(p,R_2)}$. Furthermore, if $\Omega(p,R_1,R_2) \neq \emptyset$ then the sum of the series defines a regular function $f : \Omega(p,R_1,R_2) \to \hh$.
\end{theorem}

\begin{proof} 
The Laurent series converges at all points $q$ with $\tau(q,p)>R_1$ and $\sigma(q,p) <R_2$ by the root test, since
$$\limsup_{n \to +\infty}  |(q-p)^{*n} a_n|^{1/n} = \frac{\sigma(q,p)}{R_2},$$
$$\limsup_{m \to +\infty}  |(q-p)^{-*m} a_{-m}|^{1/m} = \frac{R_1}{\tau(q,p)}.$$
For the same reason, the series does not converge at any point $q$ such that $\tau(q,p)<R_1$ or $\sigma(q,p)>R_2$. 
Now let us prove that the convergence is uniform on compact sets. Observe that
$$\Sigma(p,R_1,R_2) = \bigcup_{R_1<r_1<r_2<R_2} \overline{\Sigma(p,r_1,r_2)}$$
where $\overline{\Sigma(p,r_1,r_2)}$ is compact: it is closed by construction and it is bounded because $r_2<+\infty$. For all $q \in \overline{\Sigma(p,r_1,r_2)}$, we already know that $\sigma(q,p) \leq r_2$ and $\tau(q,p) \geq r_1$, so that $\tau(q,p)^{-m} \leq r_1^{-m}$ for all $m>0$. Hence in $\overline{\Sigma(p,r_1,r_2)}$ the function series $\sum_{n \in \nn} (q-p)^{*n} a_n$ is dominated by the convergent number series $\sum_{n \in \nn} r_2^n |a_n|$, while $\sum_{m>0} (q-p)^{-*m} a_{-m}$ is dominated by the convergent number series $\sum_{m>0} r_1^{-m} |a_{-m}|$.

Finally, let us prove the regularity of $f$ in $\Omega(p,R_1,R_2)$ when the latter is not empty. Each addend $(q-p)^{*n} a_n$ of the series defines a regular function on $\Omega(p,R_1,R_2)$. Since the convergence is uniform on compact sets, we easily deduce that the sum of the series is regular in $\Omega(p,R_1,R_2)$, as desired.
\end{proof}

Conversely, an expansion property holds. Let us begin by observing that, if $p \in \rr+I\rr$ and if $A_I (p,R_1,R_2)$ denotes the annulus $\{z \in \rr+I\rr: R_1<|z-p|<R_2\}$ then
$$\Sigma(p,R_1,R_2) = A_I(p,R_1,R_2) \cup \Omega(p,R_1,R_2).$$
Moreover, let us note what follows.

\begin{remark}
Let $p \in \hh$, let $R_1,R_2 \in \rr$ be such that $0\leq R_1<R_2\leq +\infty$ and let $I \in \s$ be such that $p \in \rr+I\rr$. If $U = A_I(p,R_1,R_2) \cap A_I(\bar p,R_1,R_2)$ denotes the largest open subset of $A_I(p,R_1,R_2)$ that is symmetric with respect to $\rr$, then $\Omega(p,R_1,R_2)$ is the symmetric completion $\widetilde{U}$ of $U$. Since $U$ has at most $2$ connected components, if $\Omega(p,R_1,R_2)\neq \emptyset$ then either $\Omega(p,R_1,R_2)$ is a symmetric slice domain or it is the disjoint union of two symmetric slice domains.
\end{remark}

We are now ready for the announced result.

\begin{theorem}
Let $f$ be a regular function on a domain $\Omega \subseteq \hh$ and let $p \in \hh$. There exists a sequence $\{a_n\}_{n \in \zz}$ in $\hh$ such that, for all $0\leq R_1<R_2\leq +\infty$ with $\Sigma(p,R_1,R_2)\subseteq \Omega$,
\begin{equation}\label{laurentformula}
f(q) = \sum_{n \in \zz} (q-p)^{*n} a_n
\end{equation}
in $\Sigma(p,R_1,R_2)$.
\end{theorem}

\begin{proof}
Suppose $p \in \rr+I\rr$: if $\Omega \supseteq \Sigma(p,R_1,R_2)$ then $\Omega_I \supseteq A_I(p,R_1,R_2)$. Chosen $J \in \s$ with $J\perp I$, let $F,G : \Omega_I \to \rr+I\rr$ be holomorphic functions such that $f_I = F+GJ$. If
$$F(z)=\sum_{n \in \zz} (z-p)^n \alpha_n,\ G(z)=\sum_{n \in \zz} (z-p)^n \beta_n$$
are the Laurent series expansions of $F$ and $G$ in $A_I(p,R_1,R_2)$ then, setting $a_n = \alpha_n+\beta_nJ$, we get
$$f_I(z)=\sum_{n \in \zz} (z-p)^n a_n$$ 
for all $z \in A_I(p,R_1,R_2)$. Since $\sum_{n \in \zz} (z-p)^n a_n$ converges in $A_I(p,R_1,R_2)$, the series in equation (\ref{laurentformula}) converges in $\Sigma(p,R_1,R_2)$. If $\Omega(p,R_1,R_2) \neq \emptyset$ then the series in equation (\ref{laurentformula}) defines a regular function $g: \Omega(p,R_1,R_2) \to \hh.$
Since $f_I \equiv g_I$ in $\Omega(p,R_1,R_2) \cap (\rr+I\rr)$, applying the Identity Principle \ref{identity} to each connected component allows us to conclude that $f$ and $g$ coincide in $\Omega(p,R_1,R_2)$ (hence in $\Sigma(p,R_1,R_2)$, as desired).
\end{proof}

It is possible to compute the coefficients of Laurent expansions by means of an integral formula. If ${\gamma_I} : [0,1] \to \rr+I\rr$ is a rectifiable curve whose support lies in a complex line $\rr+I\rr$ for some $I \in \s$, if $f$ is a continuous function on a neighborhood $\Gamma_I$ of ${\gamma_I}$ in $\rr+I\rr$, if $J \in \s$ is such that $J \perp I$  and if $F,G$ are continuous functions such that $f = F+GJ$ in $\Gamma_I$, then we denote
$$\int_{\gamma_I} ds f(s) = \int_{\gamma_I}ds F(s) +  \int_{\gamma_I}dsG(s) J.$$

\begin{proposition}
If $p \in \rr+I\rr \subset \hh$, if $0\leq R_1<R_2\leq +\infty$ and if 
$
f(q) = \sum_{n \in \zz} (q-p)^{*n} a_n
$ 
in $\Sigma(p,R_1,R_2)$, then for each $n \in \zz$
\begin{equation}
a_n = \frac{1}{2 \pi I} \int_{\gamma_I} ds (s-p)^{-n-1} f(s)
\end{equation}
where $\gamma_I(s) = p + R e^{2\pi I s}$ for any $R$ such that $R_1<R<R_2$.
\end{proposition}

\begin{proof}
Let $J \perp I$ and let $F,G$ be holomorphic functions such that $f_I=F+GJ$. If $\alpha_n,\beta_n \in \rr+I\rr$ are such that $a_n = \alpha_n + \beta_n J$ for all $n \in \zz$ then the fact that
$f_I(z) = \sum_{n \in \zz} (z-p)^{n} a_n$
for all $z$ belonging to the annulus $A_I(p,R_1,R_2)$ implies that
$$F(z) = \sum_{n \in \zz} (z-p)^{n} \alpha_n,\ G(z) = \sum_{n \in \zz} (z-p)^{n} \beta_n$$
for all $z \in A_I(p,R_1,R_2)$.
The thesis follows observing that
$$\alpha_n = \frac{1}{2 \pi I} \int_{\gamma_I} \frac{ds}{(s-p)^{n+1}} F(s)$$
$$\beta_n = \frac{1}{2 \pi I} \int_{\gamma_I} \frac{ds}{(s-p)^{n+1}} G(s).$$\qedhere
\end{proof}


\section{Classification of singularities}\label{sec:classification}

We now classify the singularities of regular functions. A first classification has been introduced in \cite{poli}, but the results of Section \ref{sec:laurent} on Laurent series allow us to change our approach.

\begin{definition}
Let $f$ be a regular function on a symmetric slice domain $\Omega$. We say that a point $p \in \hh$ is a \emph{singularity} for $f$ if there exists $R>0$ such that $\Sigma(p,0,R) \subseteq \Omega$.
\end{definition}

In other words, $p$ is a singularity for $f$ if and only if the Laurent expansion $f(q) = \sum_{n \in \zz} (q-p)^{*n}a_n$ of $f$ at $p$ has $0$ as its inner radius of convergence and it has a positive outer radius of convergence. 

\begin{remark}
For all $R>0$ and for all $p = x+Iy \in \hh$, if $\Sigma(p,R)$ denotes the $\sigma$-ball defined in Section \ref{sec:laurent} then
\begin{equation}
\Sigma(p,0,R) = \Sigma(p,R) \setminus [(x+y\s) \setminus \{\bar p\}] = [\Sigma(p,R)\setminus (x+y\s)] \cup \{\bar p\}.
\end{equation}
In other words, $\Sigma(p,0,R)$ is the union of the punctured disc $\Delta_I(p,R) \setminus \{p\}$ and of the symmetric completion $\widetilde{U}$ of $U = [\Delta_I(p,R) \setminus \{p\}] \cap [\Delta_I(\bar p,R) \setminus \{\bar p\}]$.
\end{remark}

We classify the singularities with the next Definition.

\begin{definition}
Let $p$ be a singularity for $f$ and consider the expansion
\begin{equation}\label{classification}
f(q) = \sum_{n \in \zz} (q-p)^{*n}a_n.
\end{equation}
We say that $p$ is a \emph{pole} for $f$ if there exists an $m\geq0$ such that $a_{-k} = 0$ for all $k>m$; the minimum such $m$ is called the \emph{order} of the pole and denoted by $ord_f(p)$. If $p$ is not a pole then we call it an \emph{essential singularity} for $f$ and set $ord_f(p) = +\infty$.
Finally, we say that $p$ is a \emph{removable singularity} if $f$ extends to a neighborhood of $p$ as a regular function.
\end{definition}

We point out that, when $0<R<|Im(p)|$, the set $\Sigma(p,0,R)$ reduces to the punctured disc $\Delta_I(p,R) \setminus \{p\}$ (where $I \in \s$ is chosen so that $p \in \rr+I\rr$). In this case Equation (\ref{classification}) reduces to 
$$f_I(z)=\sum_{n \in \zz} (z-p)^n a_n.$$
If $p$ is a pole for $f_I$, i.e. if there exists $n \in \nn$ such that 
$f_I(z) = (z-p)^{-n}g_I(z)$
for some holomorphic $g_I:\Delta_I(p,R) \to \hh$, it is natural to call order of $f_I$ at $p$ the least such $n$; that is, the maximum between the order of $F$ and the order of $G$ at $p$, if $f_I$ splits as $f_I=F+GJ$ according to Lemma \ref{splitting}. On the other hand, $p$ is an essential singularity for $f_I$, i.e. it is not a pole for $f_I$, if and only if it is essential for $F$ or $G$. Finally, $p$ is a removable singularity for $f_I$, i.e. $f_I$ extends as a holomorphic function to a neighborhood of $p$ in $\rr+I\rr$, if and only if $p$ is removable for both $F$ and $G$. We observe what follows.

\begin{remark}
A point $p \in \rr+I\rr \subset \hh$ is a singularity for a regular function $f$ on a symmetric slice domain if and only if it is an isolated singularity for the holomorphic function $f_I$. The function $f$ has a pole of order $n$ at $p$ if and only if $f_I$ does. Furthermore, $p$ is an essential singularity for $f$ if and only if it is essential for $f_I$.
\end{remark}

The same equivalence does not hold for removable singularities: $p$ is removable for $f_I$ if and only if it is a pole of order $0$ for $f_I$ if and only if it is a pole of order $0$ for $f$; but such a pole is not necessarily removable for $f$, as proven by the next Example.

\begin{example}
Let $I \in \s$ and let $f : \hh \setminus \s \to \hh$ be the regular function defined by 
$$f(q) = (q+I)^{-*} = (q^2+1)^{-1}(q-I).$$ 
By restricting to the complex line $\rr+I\rr$ we get $f_I(z) = \frac{1}{z+I}$ for all $z \in (\rr+I\rr) \setminus \{\pm I\}$. The point $-I$ is clearly a pole of order $1$ for $f_I$ and $f$, while $I$ is a pole of order $0$. We immediately conclude that $I$ is a removable singularity for $f_I$. The same is not true for $f$: each neighborhood $U$ of $I$ in $\hh$ includes points $q \in \s$, where $q^2+1$ vanishes while $q-I$ does not; this means that $|f|$ is unbounded in $U \setminus \s$; hence $I$ cannot be removable.
Notice that the Laurent series expansion of $f$ at $-I$ converges in 
$$\Sigma(-I,0,+\infty) = \hh \setminus (\s \setminus \{I\}) = (\hh \setminus \s) \cup \{I\}$$
while $\Omega(-I,0,+\infty) = \hh \setminus \s$ is the maximal domain on which $f$ is defined as a regular function.
\end{example}

The situation presented in this Example is quite common, as we will see in the next Section. Before proceeding, let us define the analogs of meromorphic functions.

\begin{definition}
A function $f$ is  \emph{semiregular} in a symmetric slice domain $\Omega$ if it is regular in a symmetric slice domain $\Omega'\subseteq\Omega$ such that every point of $\Omega\setminus \Omega'$ is a pole (or a removable singularity) for $f$.
\end{definition}

Notice that, if in the previous Definition we set $\mathcal{S} = \Omega\setminus \Omega'$ then, for all $I \in \s$, $\mathcal{S}_I = \mathcal{S} \cap (\rr+I\rr)$ is discrete and the restriction $f_I : \Omega_I \setminus \mathcal{S}_I \to \hh$ is a meromorphic function on $\Omega_I$. Furthermore, we observe what follows.

\begin{remark}
If $f$ is semiregular in $\Omega$ then the set $\mathcal{S}$ of its non removable poles consists of isolated real points or isolated 2-spheres of type $x+y\s$.
\end{remark}


\section{Poles and quotients}\label{sec:poles}

In this Section we study the poles and their relation to $*$-quotients. We begin by proving the semiregularity of quotients.

\begin{proposition}\label{quotientsemiregularity}
Let $f,g$ be regular functions on a symmetric slice domain $\Omega$ and consider the quotient $f^{-*}*g : \Omega \setminus Z_{f^s} \to \hh$. Each $p \in Z_{f^s}$ is a pole of order 
\begin{equation}
ord_{f^{-*}*g}(p) \leq m_{f^s}(p)
\end{equation}
for $f^{-*}*g$, where $m_{f^s}(p)$ denotes the classical multiplicity of $p$ as a zero of $f^s$. As a consequence, $f^{-*}*g$ is semiregular in $\Omega$.
\end{proposition}

\begin{proof}
Suppose $p = x+yI$ and notice that it suffices to prove the thesis for the restriction $f_I$. If $m_{f^s}(p) = n$ then, as a consequence of Proposition \ref{conjugatezeros}, $m_{f^s}(\bar p) = n$; there exists a holomorphic function $h_I$ with $h_I(p) \neq 0$ such that for all $z \in \rr+I\rr$
$$f_I^s(z) = (z-p)^n(z-\bar p)^n h_I(z)=  \left[(z-x)^2+y^2\right]^n h_I(z).$$ 
Since $f^s$ is a series with real coefficients, $h_I(z)$ must have real coefficients, too. Hence $f_I^s(z) = h_I(z) \left[(z-x)^2+y^2\right]^n = h_I(z) (z-\bar p)^n (z-p)^n$ and 
$$(f^{-*}*g)_I(z) = f_I^{-s}(z) (f^c *g)_I(z) = (z-p)^{-n}(z-\bar p)^{-n} h_I(z)^{-1}(f^c *g)_I(z)$$ 
where $(z-\bar p)^{-n} h_I(z)^{-1}(f^c *g)_I(z)$ is holomorphic in a neighborhood of $p$ in $\rr+I\rr$.
\end{proof}

Proposition \ref{quotientsemiregularity} is a direct extension of the corresponding result for $\Omega = B(0,R)$ proven in \cite{poli}. We could also extend directly the converse result, proving that a function that is semiregular in $\Omega$ can be expressed as a $*$-quotient in each $\Omega_0$ that is relatively compact in $\Omega$. However, in \cite{poli} we used this approach because we were obliged to work with Euclidean balls $B(0,R)$. Working with symmetric slice domains introduces an advantageous novelty: it is possible to treat singularities locally, in the sense of the next Lemma.

\begin{lemma}
Let $f$ be a semiregular function on a symmetric slice domain $\Omega$. If $p = x+Iy \in \Omega$ is a non removable singularity for $f$, then there exists a symmetric slice domain $U$ with $p \in U \subseteq \Omega$ such that $f$ is regular in $U \setminus (x+y\s)$.
\end{lemma}

\begin{proof}
Let $r \in \Omega \cap \rr$ be a real point such that $f$ is regular in a neighborhood of $r$ and let $\gamma : [0,1] \to \Omega_I$ be a plane curve connecting $p$ to $r$. Since the set $\mathcal{S}_I$ of non removable poles in $\Omega_I$ is discrete, we may suppose $\gamma$ not to pass through any other non removable pole of $f_I$. The image of $\gamma$ is compact, hence there exists $\varepsilon>0$ such that $\Gamma = \{z \in \rr+I\rr: \exists t \in [0,1] s.t. |z-\gamma(t)|<\varepsilon\}$ is contained in $\Omega_I$ and $f_I$ is holomorphic in $\Gamma \setminus \{p\}$. Since the set of non removable poles of $f$ is symmetric, if $U$ is the symmetric completion $\widetilde{\Gamma}$ of $\Gamma$ then $f$ must be regular in $U \setminus (x+y\s)$ where $x+y\s$ is the sphere through $p$. By construction, $U$ is a symmetric slice domain with $p \in U \subseteq \Omega$.
\end{proof}

We now prove that all semiregular functions can be locally expressed as quotients of regular functions.

\begin{theorem}\label{polefactorization1}
Let $f$ be a semiregular function on a symmetric slice domain $\Omega$. Choose $p = x+yI \in \Omega$ and set $m = ord_f(p), n = ord_f(\bar p)$: without loss of generality $m \leq n$. There exist a neighborhood $U$ of $p$ in $\Omega$ that is a symmetric slice domain and a (unique) regular function $g : U \to \hh$ such that
\begin{eqnarray}
f(q) &=& [(q-p)^{*m}*(q-\bar p)^{*n}]^{-*}* g(q) = \\ 
 &=& \left[(q-x)^2+y^2\right]^{-n} (q-p)^{*(n-m)}* g(q) \nonumber
\end{eqnarray}
in $U \setminus (x+y\s)$. Moreover, if $n>0$ then $g(p) \neq 0,g(\bar p) \neq 0$.
\end{theorem}

\begin{proof}
Choose $U$ as in the previous Lemma. The restriction $f_I$ is meromorphic in $U_I$, it is holomorphic in $U_I \setminus \{p,\bar p\}$ and it has a pole of order $m$ at $p$ and a pole order of $n$ at $\bar p$. Hence there exists a holomorphic $g_I : U_I \to \hh$ with $$f_I(z) = \frac{1}{(z-p)^m(z-\bar p)^n}g_I(z)$$ 
for all $z \in U_I \setminus \{p,\bar p\}$. Let $g = \ext(g_I)$ and consider the function 
$$h(q) = [(q-p)^{*m}*(q-\bar p)^{*n}]^{-*}* g(q).$$ 
$h$ is regular on its domain of definition, which is $U \setminus (x+y\s)$. Furthermore, $h_I (z) = \frac{1}{(z-p)^m(z-\bar p)^n}g_I(z) =f_I(z)$ for all $z \in U_I \setminus \{p,\bar p\}$. The Identity Principle \ref{identity} allows us to conclude that $f(q) = h(q)$ for all $q \in U \setminus (x+y\s)$.
The second equality is proven observing that 
$$[(q-p)^{*m}*(q-\bar p)^{*n}]^{c} = (q-p)^{*n}*(q-\bar p)^{*m} =$$
$$= [(q-x)^2+y^2]^{m}  (q-p)^{*(n-m)}$$ 
and $[(q-p)^{*m}*(q-\bar p)^{*n}]^{-s} = [(q-x)^2+y^2]^{-m-n}$, so that 
$$[(q-p)^{*m}*(q-\bar p)^{*n}]^{-*}= \left[(q-x)^2+y^2\right]^{-n} (q-p)^{*(n-m)}.$$
\end{proof}

We derive what follows.

\begin{proposition}
The set of semiregular functions on a symmetric slice domain $\Omega$ is a division ring with respect to $+,*$.
\end{proposition}

\begin{proof}
Let $f,g$ be two semiregular functions on $\Omega$ and let $\mathcal{S}_f, \mathcal{S}_g$ be the sets of non removable singularities of $f,g$ (respectively). The sum $f+g$ and the $*$-product $f*g$ are defined (and regular) on the largest symmetric slice domain in which both $f$ and $g$ are regular, $\Omega \setminus (\mathcal{S}_f \cup \mathcal{S}_g)$.
Moreover, if $f \not \equiv 0$ then $f^{-*}$ is defined on $\Omega \setminus (\mathcal{S}_f \cup Z_{f^s})$. The functions $f+g, f*g$ are semiregular in $\Omega$ since every $p \in \mathcal{S}_f \cup \mathcal{S}_g$ is a pole for $f+g$ and $f*g$. Indeed, by Theorem \ref{polefactorization1}, there exists a neighborhood $U$ of $p$ in $\Omega$ where $f$ and $g$ can be expressed as quotients of regular functions. By Theorem \ref{quotients}, $f+g$ and $f*g$ are quotients of regular functions on $U$, too. In particular, by Proposition \ref{quotientsemiregularity}, $f+g$ and $f*g$ are semiregular in $U$. Finally, if $f \not \equiv 0$ then $f^{-*}$ is semiregular in $\Omega$ since every $p \in \mathcal{S}_f \cup Z_{f^s}$ is a pole for $f^{-*}$ by a similar reasoning.
\end{proof}

We can now state and prove the following consequence of Theorem \ref{polefactorization1}.

\begin{corollary}\label{polefactorization}
Let $f$ be a semiregular function on a symmetric slice domain $\Omega$. Choose $p = x+yI \in \Omega$, let $m = ord_f(p), n = ord_f(\bar p)$ and suppose $m \leq n$. There exists a unique semiregular function $g$ on $\Omega$, without poles in $x+y\s$, such that
\begin{eqnarray}
f(q) &=& [(q-p)^{*m}*(q-\bar p)^{*n}]^{-*}* g(q) = \\ 
 &=& \left[(q-x)^2+y^2\right]^{-n} (q-p)^{*(n-m)}* g(q) \nonumber
\end{eqnarray}
Furthermore, if $n>0$ then $g(p) \neq 0,g(\bar p) \neq 0$.
\end{corollary}

\begin{proof}
By the previous Proposition, setting $g(q) = (q-p)^{*m}*(q-\bar p)^{*n}*f(q)$ defines a function which is semiregular in $\Omega$. The thesis follows by Theorem \ref{polefactorization1}.
\end{proof}

The next result explains the distribution of the poles of semiregular functions in each 2-sphere $x+y\s$.

\begin{theorem}[Structure of the poles]
Let $\Omega$ be a symmetric slice domain and let $f$ be semiregular in $\Omega$. In each sphere $x+y\s \subset \Omega$, all the poles have the same order $ord_f$ with the possible exception of one, which must have lesser order.
\end{theorem}

\begin{proof}
Choose a sphere $x+y\s \subset \Omega$ for which there exists $I \in \s$ such that $p = x+yI$ and $\bar p = x-yI$ have orders $m$ and $n$ with $m>0$ or $n>0$. Without loss of generality, $m \leq n$. By Corollary \ref{polefactorization}, there exists a semiregular function $g$ on $\Omega$ which is regular in a neighborhood $U$ of $x+y\s$ such that 
$$f(q) = [(q-x)^2+y^2]^{-n} (q-p)^{*(n-m)}* g(q)$$ 
and $g(p),g(\bar p) \neq 0$. If $\tilde f(q) = (q-p)^{*(n-m)}* g(q)$ then 
$$f(q) = [(q-x)^2+y^2]^{-n} \tilde f(q),$$ 
$$f_J(z) = [z-(x+yJ)]^{-n} [z-(x-yJ)]^{-n} \tilde f_J(z)$$ 
for all $J \in \s$. Now:
\begin{enumerate}
\item If $m<n$ then $\tilde f(x+yI) = 0$ and $\tilde f (x+yJ) \neq 0$ for all $J \in \s \setminus \{I\}$. The last equality allows us to conclude $ord_f(x+yJ) = n$ for all $J \in \s \setminus \{I\}$. Since we know by hypothesis that $ord_f(x+yI)= ord_f(p) = m<n$, the thesis holds.
\item If $m = n$ then $\tilde f(x+yI) \neq 0$. If $\tilde f$ does not have zeros in $x +y \s$ then we conclude $ord_f(x+yJ) = n$ for all $J \in \s$. If, on the contrary, $\tilde f(x + yK)=0$ for some $K \in \s$ then we can factor $z-(x+yK)$ out of $\tilde f_K(z)$ and conclude that $ord_f(x + yK)<n$ while $ord_f(x + yJ)=n$ for all $J \in \s \setminus \{K\}$, as desired. Notice that, by construction, $\tilde f$ cannot have more than one zero in $x+y\s$.
\qedhere
\end{enumerate}
\end{proof}

The previous result inspired the following.

\begin{theorem}
Let $f$ be a semiregular function on a symmetric slice domain $\Omega$, suppose $f \not \equiv 0$ and let $x+y\s \subset \Omega$. There exist $m \in \zz, n \in \nn$, $p_1,...,p_n \in x+y\s$ with $p_i \neq \bar p_{i+1}$ for all $i \in \{1,\ldots, n-1\}$, so that
\begin{equation}
f(q) = [(q-x)^2+y^2]^m (q-p_1)*(q-p_2)*...*(q-p_n)*g(q)
\end{equation}
for some semiregular function $g$ on $\Omega$ which does not have poles nor zeros in $x+y\s$. If $m\leq0$ then we say that $f$ has \emph{spherical order} $-2m$ at $x+y\s$. Whenever $n>0$, we say that $f$ has \emph{isolated multiplicity} $n$ at $p_1$.
\end{theorem}

\begin{proof}
The thesis is an immediate consequence of Theorem \ref{factorization} and of Corollary \ref{polefactorization}.
\end{proof}

The spherical order and isolated multiplicity of $f$ are related to the order $ord_f$ in the following way.

\begin{remark}
Let $f$ be a semiregular function on $\Omega$ which is not regular at $p = x+yI \in \Omega$. Then $f$ has spherical order $2 \max \{ord_f(p),ord_f(\bar p) \}$ at $x+y\s$. If moreover $ord_f(p)>ord_f(\bar p)$, then $f$ has isolated multiplicity $n \geq ord_f(p)-ord_f(\bar p)$ at $\bar p$.
\end{remark}


\section{Casorati-Weierstrass Theorem}\label{sec:casorati}

In this Section we study essential singularities, proving a quaternionic version of the Casorati-Weierstrass Theorem. We begin by showing that the $*$-quotient $f^{-*}*g(q)$ is nicely related to the pointwise quotient $f(q)^{-1}g(q)$ (generalizing \cite{zerosopen,open,poli}).

\begin{proposition}\label{reciprocalformula}
Let $f,g$ be regular functions on a symmetric slice domain $\Omega$. Then, for all $q \in \Omega \setminus Z_{f^s}$,
\begin{equation}
f^{-*}*g(q) = f(T_f(q))^{-1} \cdot g(T_f(q)),
\end{equation}
where $T_f : \Omega \setminus Z_{f^s} \to \Omega \setminus Z_{f^s}$ is defined by $T_f(q) = f^c(q)^{-1} q f^c(q)$. Furthermore, $T_f$ and $T_{f^c}$ are mutual inverses so that $T_f$ is a diffeomorphism.
\end{proposition}

\begin{proof}
If $f^s(q) \neq 0$ then $f^c(q)\neq 0$. Hence $T_f$ is well defined on $\Omega \setminus Z_{f^s}$, thanks to Proposition \ref{conjugatezeros}. Recalling Proposition \ref{formprod}, we compute:
$$f^{-*}*g(q) = f^s(q)^{-1} f^c*g(q) = \left[ f^c*f(q) \right]^{-1} f^c(q) g(T_f(q))=$$ 
$$= \left[f^c(q) f(T_f(q))\right]^{-1} f^c(q)  g(T_f(q))= f(T_f(q))^{-1} f^c(q)^{-1} f^c(q)  g(T_f(q)) =$$
$$= f(T_f(q))^{-1}  g(T_f(q)).$$
Moreover, $T_f : \Omega \setminus Z_{f^s} \to \hh$ maps any sphere $x+y\s$ to itself. In particular, since $Z_{f^s}$ is symmetric by Proposition \ref{conjugatezeros}, $T_f(\Omega \setminus Z_{f^s}) \subseteq \Omega \setminus Z_{f^s}$.
Now, since $(f^c)^c = f$ we observe that $T_{f^c}(q) = f(q)^{-1} q f(q)$. For all $q \in \Omega \setminus Z_{f^s}$,  and setting $p = T_f(q)$ we have that
$$T_{f^c} \circ T_f(q) = T_{f^c}(p) = f(p)^{-1} p f(p) =$$ 
$$= f(p)^{-1} \left[f^c(q)^{-1} q f^c(q)\right] f(p) = \left [f^c(q)f(p)\right]^{-1} q \left [f^c(q)f(p)\right]$$ 
where
$$f^c(q)f(p) = f^c(q) f(f^c(q)^{-1}q f^c(q)) = f^c* f (q) = f^s(q).$$ 
Now, $f^s(q)$ and $q$ always lie in the same complex line $\rr+I\rr$. In particular they commute, so that
$$T_{f^c} \circ T_f(q) =  f^s(q)^{-1}  q  f^s(q) = q,$$
as desired.
\end{proof}

We are now ready for the announced result.

\begin{theorem}[Casorati - Weierstrass]
Let $\Omega$ be a symmetric slice domain and let $f$ be a regular function on $\Omega$. If $p$ is an essential singularity for $f$ and if $U$ is a symmetric neighborhood of $p$ in $\hh$, then $f(\Omega \cap U)$ is dense in $\hh$.
\end{theorem}

\begin{proof}
Suppose that for some symmetric neighborhood $U$ of $p$ in $\hh$ there existed a $v \in \hh$ and an $\varepsilon >0$ such that $f(\Omega \cap U) \cap B(v, \varepsilon) = \emptyset$. Setting $h(q) = f(q) - v$ and $g(q) = h^{-*}$ would then define a function $g$, semiregular in $\Omega$, with
$$|g(q)|= \frac{1}{|f(T_h(q)) - v|} \leq \sup_{w \in \Omega \cap U} \frac{1} {|f(w) - v|} \leq \frac{1}{\inf_{w \in \Omega \cap U} |f(w) - v|} \leq \frac{1}{\varepsilon}$$ 
for all $q \in \Omega \cap U$. The function $g$ would then have a removable singularity at $p$ and the function $f = g^{-*} + v$, would have a pole (or a removable singularity) at $p$. This is impossible, since we supposed $p$ to be an essential singularity for $f$.
\end{proof}

The previous Theorem was proven in the special case $p=0$ in \cite{moebius}. We notice that in this case, and in general when $p \in \rr$, the situation is completely analogous to the complex setting.

\begin{corollary}
Let $f$ be a regular function on $B(p,R) \setminus \{p\}$ with $p \in \rr, R>0$. If $f$ has an essential singularity at $p$ then for each neighborhood $U$ of $p$ in $B(p,R)$, the set $f(U \setminus \{p\})$ is dense in $\hh$.
\end{corollary}




\bibliographystyle{abbrv}
\bibliography{Singularities}


\end{document}